\documentclass[letterpaper, 10 pt, conference]{ieeeconf}   % Comment this line out
\IEEEoverridecommandlockouts
\usepackage{graphics} % for pdf, bitmapped graphics files
\usepackage{epsfig} % for postscript graphics files
\usepackage{mathptmx} % assumes new font selection scheme installed
\usepackage{times} % assumes new font selection scheme installed
\usepackage{amsmath} % assumes amsmath package installed
\usepackage{amssymb}  % assumes amsmath package installed
\usepackage{psfrag}
\title{\LARGE \bf
Positive Consensus of Directed Multi-agent Systems
}

\author{Nachuan Yang,~\IEEEmembership{Student Member,~IEEE}, Yonghua Yin and Jinrong Liu$^{\dag}$% <-this % stops a space
\thanks{Nachuan Yang is with the Department of Mathematics, Faculty of Science,
        The University of Hong Kong, Pokfulam Rd, Hong Kong
        {\tt\small yangnachuan@connect.hku.hk}}%
\thanks{Yonghua Yin is with the Department of Electrical and Electronic Engineering, Imperial College London, UK, SW7 2BT
        {\tt\small y.yin14@imperial.ac.uk}}
\thanks{Jinrong Liu is with the Department of Mechanical Engineering, Faculty of Engineering,
        The University of Hong Kong, Pokfulam Rd, Hong Kong
        {\tt\small liujinrjason@connect.hku.hk}}%
\thanks{$^{\dag}$ indicates corresponding author}
}

\begin{document}

\maketitle
\thispagestyle{empty}
\pagestyle{empty}

%%%%%%%%%%%%%%%%%%%%%%%%%%%%%%%%%%%%%%%%%%%%%%%%%%%%%%%%%%%%%%%%%%%%%%%%%%%%%%%%
\begin{abstract}

This paper addresses the problem of positive consensus of directed multi-agent systems with observer-type output-feedback protocols. More specifically, directed graph is used to model the communication topology of the multi-agent system and linear matrix inequalities (LMIs) are used in the consensus analysis in this paper. Using positive systems theory and graph theory, a convex programming algorithm is developed to design appropriate protocols such that the multi-agent system is able to reach consensus with its state trajectory always remaining in the non-negative orthant. Finally, numerical simulations are given to illustrate the effectiveness of the derived theoretical results.

\end{abstract}

%%%%%%%%%%%%%%%%%%%%%%%%%%%%%%%%%%%%%%%%%%%%%%%%%%%%%%%%%%%%%%%%%%%%%%%%%%%%%%%%
\section{INTRODUCTION}
%1. Consensus of multi-agent systems.
Recently, the research on multi-agent systems, especially the consensus issue has received much attention from various scientific and engineering areas for its important applications in sensor networks, automatic vehicles and modern robotics, to name just a few \cite{Ahmed2014Asynchronous, Fax2004Information, Mu2016Design}. Theoretically speaking, the main concern on this issue is to design effective control protocols so that the multiple agents in the overall system are able to cooperatively and coordinately attain some common goals without centralized controllers. Most existing studies on the consensus issue of multi-agent systems assume that full state information of agents is known. Based on this assumption \cite{Valcher2014On, Valcher2017On}, many algorithms have been developed to design static consensus protocols. Recently, dynamic output-feedback protocols have been broadly used
to solve the consensus problem of multi-agent systems \cite{wang2017observer, Li2010Consensus}.
Hence, how to design such kind of control protocols has become an important issue nowadays.

%2. Positive systems
Positive systems have the special property that, the states and outputs of a positive system are always non-negative if its initializations and inputs are non-negative. The applications of positive systems can be very broad, including industrial engineering, systems biology, and biomedicine \cite{Caccetta2004A, farina2011positive, Jong2003Hybrid}. Among quantities of research literature on positive systems, special attention has been devoted to the reachability and realization of such kind of systems. For instance, sufficient and necessary conditions on positive realizability have been concluded via convex analysis in \cite{farina1996On}. The synthesis and analysis on positive dynamics are investigated using linear matrix inequality (LMI) method and new results are concluded in \cite{Ebihara2014LMI}. In recent years, positive systems theory is also applied in the study of nodal networks, time delay system and edge-consensus problem, and many useful results have been derived on these problems \cite{ebihara2017steady, Han2017Discrete, Su2018Positive}. For multi-agent systems, positive systems are commonly used to model the dynamics of the agents. A classical example is the multi-agent system with integrators as agents where the multiple agents are regarded as positive systems \cite{Fax2004Information}. There are also lots of other examples where positive multi-agent systems are involved \cite{Knorn2011A, Pahuja2013A}. In real applications, values involved with practical systems are usually intrinsically non-negative, so the positivity should be guaranteed when analyzing the consensus issue of such kinds of multi-agent systems \cite{Knorn2011A}. For these reasons, we are motivated to investigate the positive consensus of multi-agent systems.

%3. Positive consensus
Although many breakthroughs on positive consensus of multi-agent systems have been made in the past few years, the complete solution to this challenging problem is still under investigation. In a common way, the general consensus problem can be transformed to a simultaneous stabilization problem, but this transformation cannot be directly applied to positive multi-agent systems because the positivity of the overall system cannot always be guaranteed \cite{Valcher2014On}. Recently, many new results have been obtained on this problem. In \cite{Valcher2014On}, a single-input single-output positive state-space model is used to describe the agents of multi-agent systems, and some conditions on positive consensus are concluded. This problem is further discussed in \cite{Li2010Consensus, Valcher2016New, Valcher2017On} where undirected multi-agent systems are considered. In the recent work \cite{HanObserver}, the consensus of positive multi-agent systems with strongly connected and directed communication topology is studied. These works provide many useful results to solve the positive consensus problem of multi-agent systems. However, the positive consensus of multi-agent systems with general directed communication topology is still an open question. This motivates our work in this paper. Compared with the existing work \cite{HanObserver} where the multi-agent systems are assumed to be directed and strongly connected, we investigate the positive consensus issue of multi-agent systems in a more general case, where the communication topology is directed and only assumed to have a spanning tree.
%%4. Contributions
%The contributions of this paper are concluded in the following points. To design observer-type dynamic output-feedback control protocols such that the positive consensus of multi-agent systems can be guaranteed, we derive some new theoretical results and develop an effective convex programming algorithm. Some tools in graph theory, such as the Laplacian matrix, also play an important role to derive the final results.

%5. Organizations
The rest of this paper is organized as follows. In Section 2, some background and preliminaries on positive systems theory and graph theory are provided and the problem studied in this paper is defined. In Section 3, consensus analysis and design of positive multi-agent systems with observer-type dynamic protocols are derived and a convex programming algorithm is developed. In Section 4, numerical simulations are given to illustrate the effectiveness of the algorithm. In Section 5, the whole paper is summarized and concluded.

%%%%%%%%%%%%%%%%%%%%%%%%%%%%%%%%%%%%%%%%%%%%%%%%%%%%%%%%%%%%%%%%%%%%%%%%%%%%%%%%
\section{Notations and Preliminaries}

\subsection{Notations}

In this paper, capital letters such as $A$ are used to denote matrices and lowercase letters such as $v$ represent scalars, or vectors if stated that $v \in \mathbb{R}^{m}$. For scalar $v \in \mathbb{C}$, the notation ${\rm Re}(v)$ means the real part of $v$. The notation $A \in \mathbb{R}^{m\times n}$ means that, all entries of matrix $A$ are real and matrix $A$ has $m$ rows and $n$ columns. Matrices in this paper are assumed to have compatible dimensions if not explicitly stated. $I_{n}$ denotes the $n\times n$ identity matrix and $I$ denotes the identity matrix with appropriate dimension. $\mathbf{1}_{n}$ denotes the $n$-dimensional vector whose all entries are equal to one. The superscript \text{T} represents the transpose of a matrix. The superscript \text{*} represents the conjugate transpose of a matrix. For matrix $A \in \mathbb{R}^{m\times n}$, $[A]_{ij}$ denotes the
element located at the $i$-th row and the $j$-th column. The notation $A\succeq 0$ (respectively, $A\succ0$) means that for all $i$ and $j$, $[A]_{ij} \succeq 0$ (respectively, $[A]_{ij} \succ 0$). The notation $A \succeq B$ (respectively,
$A \succ B$) means that the matrix $A-B \succeq 0$ (respectively,
$A-B \succ 0$). The notation $A \geq 0$ (respectively, $A>0$) means that A is positive semidefinite (respectively, positive definite). The notation $A \geq B$ (respectively, $A>B$) means that the matrix $A-B$ is positive semi-definite (respectively, positive definite). The notation $A \otimes B$ denotes the Kronecker product of matrices $A$ and $B$. Matrix $A\in \mathbb{R}^{n\times n}$ is called Metzler if all its off-diagonal elements are non-negative, i.e., $[A]_{ij} \succeq 0$ whenever $i \neq j$, which is denoted by $A \in \mathbb{M}^{n}$.
Matrix $A\in \mathbb{R}^{n\times n}$ is called Hurwitz if all its eigenvalues have strictly negative real part, i.e., ${\rm Re}(\lambda_{i}) \prec 0$ for each eigenvalue $\lambda_{i}$, which is denoted by $A \in \mathbb{H}^{r}$. The notation $\alpha(A)$ means the spectral abscissa of matrix $A$.

\subsection{Positive Systems Theory}
Consider a continuous-time linear system:
\begin{equation}\label{LTIsystem}
\begin{cases}
\dot{x}(t)=Ax(t)+Bu(t)  \\
y(t)=Cx(t)
\end{cases}
\end{equation}
where $x(t)\in \mathbb{R}^{r}$, $u(t)\in \mathbb{R}^{m}$, and $y(t)\in \mathbb{R}^{p}$ denote
the system state, control input and output respectively. $A$, $B$ and $C$ are system matrices with compatible dimensions.
In order to analyze the positive consensus of multi-agent systems, some useful results are listed as follows \cite{farina2011positive}:
\newtheorem{myDef}{Definition}
\begin{myDef}\label{definition1}
System (\ref{LTIsystem}) is a continuous-time
linear positive system if for all initial value $x(0)\succeq 0$ and input $u(t)\succeq 0$, then the state trajectory $x(t)\succeq 0$,  and the output $y(t)\succeq 0$ for all $t \succeq 0$.
\end{myDef}
\newtheorem{lemma}{Lemma}
\begin{lemma}\label{Plemma1}
System (\ref{LTIsystem}) is positive if and
only if matrix $A$ is Metzler, matrices $B$ and $C$ are non-negative, i.e., $A \in \mathbb{M}^{r}$, $B \succeq 0$, and $C \succeq 0$.
\end{lemma}
\vskip 0.2 cm
\begin{lemma}\label{Plemma2}
If system (\ref{LTIsystem}) is a continuous-time linear positive system, then it is asymptotically stable if and only if there exists a diagonal matrix $D>0$ such that
\begin{equation*}
A^{\mathrm{T}}D+DA<0 ~~ or ~~ DA^{\mathrm{T}}+AD<0.
\end{equation*}
\end{lemma}

\subsection{Graph Theory}
Graph is used to describe the communication topology of multi-agent systems. Mathematically speaking, graph is a structure composed of vertices where some of them are connected by edges. If all edges in a graph have no orientation, it is called undirected graph. Otherwise, it is called directed graph.  A path in a graph is a sequence of end-to-end (directed) edges. Without loss of generality, directed graph is used to model the communication topology of multi-agent systems in this paper. A graph can be described by an ordered set $\mathcal{G } =(\mathcal{V}$, $\mathcal{E} )$ which consists of a finite vertex set $\mathcal{V}=\{ $$v_{1}$, $v_2$, $\ldots$ , $v_n \}$  and an edge set $\mathcal{E} \subset \mathcal{V} \times \mathcal{V}$. For convenience, we also define a number set $\mathcal{I} := \{1, 2, $\ldots$, n \}$. A directed graph is said to have a spanning tree if there is a node such that there exists a path from it to any other nodes in the graph. For the purpose of consensus, $\mathcal{G }$ is assumed to have a spanning tree in this paper. Degree matrix for graph $\mathcal{G }$ is defined as a diagonal matrix $\mathcal{D}$ where $[\mathcal{D}]_{ii}$ is equal to the indegree of $v_{i}$, i.e., the number of incoming edges at $v_{i}$. Adjacency matrix for directed graph $\mathcal{G }$ is defined as an $n \times n$ matrix $\mathcal{A}$ such that $[\mathcal{A}]_{ij} = 1$ if there is a directed edge from $v_{j}$ to $v_{i}$, i.e., $(v_{j}, v_{i}) \in \mathcal{E}$; otherwise, $[\mathcal{A}]_{ij} = 0$. In this paper, we also assume that graph $\mathcal{G}$ has no self-loop, i.e., $[\mathcal{A}]_{ii} = 0$, $i \in \mathcal{I}$. If $(v_{j}$, $v_{i})$ $\in \mathcal{E}$, then $v_{j}$ is called the neighbour of $v_{i}$. The set of all neighbours of $v_{i}$ is denoted by $\mathcal{N}_{i}=\{v_{j}\in \mathcal{V}: (v_j,v_i)\in \mathcal{E} \}$.
The graph $\mathcal{G}$ can also be described by its Laplacian matrix $\mathcal{L}$. The Laplacian matrix $\mathcal{L}$ for graph $\mathcal{G}$ is defined as $\mathcal{L}:= \mathcal{D} - \mathcal{A}$, i.e., $[\mathcal{L}]_{ii}=  \sum_{v_j \in \mathcal{N}_{i}} [\mathcal{A}]_{ij}  $ and  $[\mathcal{L}]_{ij}= -[\mathcal{A}]_{ij} $ for any $i\neq j$. Define $l_{\mathrm{max}}:=\mathrm{max}( [\mathcal{L}]_{ii})$, $i\in \mathcal{I}$. It is easy to see that each row of $\mathcal{L}$ sums up to $0$, and thus $\mathbf{1}_{n}$ is always an eigenvector of $\mathcal{L}$ corresponding to the eigenvalue $0$. As graph $\mathcal{G}$ is directed, the eigenvalues of $\mathcal{L}$ are complex numbers, which can be ordered and denoted as $0 = \lambda_{1} \prec {\rm Re}(\lambda_{2}) \preceq \ldots \preceq {\rm Re}(\lambda_{n})$.

\subsection{Problem Formulation}
Consider a multi-agent system with $n$ identical agents, where each agent has linear positive dynamics. It can be described by
\begin{equation}\label{nodes}
\begin{cases}
\dot{x}_{i}(t)=Ax_{i}(t)+Bu_{i}(t) \\
y_{i}(t)=Cx_{i}(t),~~i \in  \mathcal{I}
\end{cases}
\end{equation}
where $x_{i}(t):=[x_{i1}$, $x_{i2}$, $\ldots$ , $x_{ir}]^{\mathrm{T}}\in \mathbb{R}^{r}$ is the state, $u_{i}(t)\in \mathbb{R}^{m}$ is the control input, $y_{i}(t)\in \mathbb{R}^{p}$ is the measured output. $A \in \mathbb{R}^{r \times r}$ is a Metzler matrix,  $B  \in \mathbb{R}^{r \times m}$ and $C  \in \mathbb{R}^{p \times r}$ are non-negative matrices. Moreover, $(A,~ B,~ C)$ is assumed to be detectable and stabilizable
in this paper.

The following observer-type dynamic output-feedback protocol is used:
\begin{equation}\label{OBOF}
\begin{cases}
\dot{\hat{x}}_{i}(t)=A\hat{x}_{i}(t)+ LC\sum_{v_j\in\mathcal{N}_{i}} [\mathcal{A}]_{ij}({e}_{j}(t) - {e}_{i}(t)) + Bu_{i}(t)\\
u_{i}(t)= - K{\hat{x}}_{i}(t)
\end{cases}
\end{equation}
where $i \in  \mathcal{I}$, $\hat{x}_{i}(t)\in \mathbb{R}^{r}$ is the state, $u_{i}(t)\in \mathbb{R}^{m}$ is the control input and $e_{i}(t)$ is the feedback signal for agent $i$ which is defined as $e_{i}(t):=x_{i}(t)-\hat{x}_{i}(t)$. $K$ and $L$ are feedback gain matrices to be determined.

Using $\tilde{x}_{i}(t):=[{x}_{i}^{\mathrm{T}}(t)$, $e_{i}^{\mathrm{T}}(t)]^{\mathrm{T}}$ as the state variable, the augmented system for each agent $i$ with the observer-type dynamic protocol can be described as

\begin{equation}\label{SOF}
\begin{cases}
\dot{\tilde{x}}_{i}(t)
= \tilde{A}\tilde{x}_{i}(t)+\tilde{B}\tilde{u_{i}}(t) \\
\tilde{y}_{i}(t) = \tilde{C}\tilde{x}_{i}(t) \\
\tilde{u}_{i}(t) = L\sum_{v_j\in\mathcal{N}_{i}} [\mathcal{A}]_{ij}(\tilde{y}_{j}(t)-\tilde{y}_{i}(t))
\end{cases}
\end{equation}
where
\begin{equation}\label{ABC}
\tilde{A}=
\begin{bmatrix}
	A - BK     &    BK \\
	0      &    A
\end{bmatrix},~~~
\tilde{B}=
\begin{bmatrix}
	0       \\
	I
\end{bmatrix},~~~
\tilde{C}=
\begin{bmatrix}
	0    & C
\end{bmatrix}.
\end{equation}

Define the state $X(t):=[\tilde{x}_{1}^{\mathrm{T}}(t)$, $\tilde{x}_{2}^{\mathrm{T}}(t)$, $\ldots$ , $\tilde{x}_{n}^{\mathrm{T}}(t)]^{\mathrm{T}} \in\mathbb{R} ^{rn}$.
Then the overall closed-loop system is represented by
\begin{equation}\label{Closedloop}
\dot{X}(t)=\mathbf{A} X(t)
\end{equation}
where $\mathbf{A}=I_{n}\otimes \tilde{A}- \mathcal{L}\otimes \tilde{B}L\tilde{C} $.

The positive consensus problem of directed multi-agent systems is studied in this paper. Based on the above descriptions, the problem to be solved is defined as follows:
\vskip 0.2 cm
\noindent
\textbf{Problem PCDMAS} (Positive Consensus of Directed Multi-agent Systems): Regarding a multi-agent system (2) with observer-type dynamic output-feedback control protocol (3), assuming that all agents have identical positive dynamics, given any non-negative initial values, design matrices $K$ and $L$ such that the consensus of the nominal dynamic system in (\ref{SOF}) and (\ref{ABC}) is achievable, i.e., $\lim_{t \rightarrow \infty}(\tilde{x}_{j}(t)-\tilde{x}_{i}(t))=0,~\forall i,j \in \mathcal{I}$, meanwhile the state of each augmented system in the overall closed-loop system keeps non-negative, i.e., $X(t) \succeq 0$ for $t \succeq 0$.
\vskip 0.2 cm
\newtheorem{remark}{Remark}
\begin{remark}
It can be observed from (\ref{SOF}) and (\ref{ABC}) that, the designs of the feedback gain matrices $L$ and $K$ are separated in the overall system so that the matrices $L$ and $K$ can be designed independently \cite{luenberger}. Moreover, as the matrices $K$ and $L$ should be designed such that $\hat{x}_{i}(t)$ converges to zero asymptotically, i.e., $\lim_{t \rightarrow \infty}\hat{x}_{i}(t)=0, ~\forall i \in \mathcal{I}$, the consensus problem of the multi-agent system in (\ref{nodes}) via the observer-type dynamic output-feedback protocol in (\ref{OBOF}) can be transformed to the consensus problem of the augmented system in (\ref{SOF}) and (\ref{ABC}).
\end{remark}

%%%%%%%%%%%%%%%%%%%%%%%%%%%%%%%%%%%%%%%%%%%%%%%%%%%%%%%%%%%%%%%%%%%%%%%%%%%%%%%%
\section{Main Results}
In this section, \textbf{Problem PCDMAS} is studied based on the results of positive systems theory and the consensus issue. Some new results on positive consensus of directed multi-agent systems are derived and a convex programming algorithm is developed to design the protocols.
\vskip 0.2 cm
\newtheorem{myTheo}{Theorem}
\begin{myTheo}
\textbf{Problem PCDMAS} is solvable if and only if all the following conditions hold:

\noindent
{\rm 1)} $BK \succeq 0$,

\noindent
{\rm 2)} $A-BK$ is Metzler and Hurwitz,

\noindent
{\rm 3)} $LC \succeq 0$,

\noindent
{\rm 4)} $A-l_{\mathrm{max}}LC$ is Metzler,

\noindent
{\rm 5)} $A-\lambda_{i}LC$ is Hurwitz, $ \forall i \in \mathcal{I} \backslash \{ 1\}  $.
\end{myTheo}
\textbf{Proof. }

(i) Positivity:  According to \textbf{\rm Lemma 1}, the overall closed-loop system (6) is positive if and only if the system matrix $\mathbf{A}$ is Metzler. By definition, the system matrix $\mathbf{A}$ can be represented as $\mathbf{A}=$
\begin{equation}\label{Amatrix1}
\begin{bmatrix}
\begin{smallmatrix}
\tilde{A}-\sum_{v_j \in \mathcal{N}_{1}} [\mathcal{A}]_{1j}\tilde{B}L\tilde{C} &  [\mathcal{A}]_{12}\tilde{B}L\tilde{C}                 & \ldots & [\mathcal{A}]_{1n}\tilde{B}L\tilde{C}  \\
[\mathcal{A}]_{21}\tilde{B}L\tilde{C}                 &\tilde{A}-\sum_{v_j \in \mathcal{N}_{2}} [\mathcal{A}]_{2j}\tilde{B}L\tilde{C}    & \ldots & [\mathcal{A}]_{2n}\tilde{B}L\tilde{C}\\
\vdots                        &\vdots                           & \ddots & \vdots\\
[\mathcal{A}]_{n1}\tilde{B}L\tilde{C}                 & [\mathcal{A}]_{n2}\tilde{B}L\tilde{C}                 & \ldots & \tilde{A}-\sum_{v_j \in \mathcal{N}_{n}} [\mathcal{A}]_{nj}\tilde{B}L\tilde{C}
\end{smallmatrix}
\end{bmatrix}
\end{equation}
where
\begin{equation}\label{Amatrix2}
\tilde{A}-\sum_{v_j \in \mathcal{N}_{i}} [\mathcal{A}]_{ij} \tilde{B}L\tilde{C}=
\begin{bmatrix}
	A - BK      &    BK \\
	0           &    A-\sum_{v_j \in \mathcal{N}_{i}} [\mathcal{A}]_{ij}LC
\end{bmatrix}
\end{equation}
and
\begin{equation}\label{aBK}
[\mathcal{A}]_{ij}\tilde{B}L\tilde{C}=
\begin{bmatrix}
        0      &    0 \\
       	0     &    [\mathcal{A}]_{ij}LC
\end{bmatrix}.
\end{equation}
It is easy to see that, $\mathbf{A}$ is Metzler if and only if  $[\mathcal{A}]_{ij}\tilde{B}L\tilde{C}$ is non-negative and $\tilde{A}-\sum_{v_j \in \mathcal{N}_{i}} [\mathcal{A}]_{ij} \tilde{B}L\tilde{C}$ is Metzler. Since $[\mathcal{A}]_{ij} \succeq 0$, we have $LC \succeq 0$ by (9). From equation (8), $\tilde{A}-\sum_{v_j \in \mathcal{N}_{i}} [\mathcal{A}]_{ij} \tilde{B}L\tilde{C} \in \mathbb{M}^{r}$ is equivalent to, $BK \succeq 0$, $A - BK \in \mathbb{M}^{r}$ and $A-\sum_{v_j \in \mathcal{N}_{i}} [\mathcal{A}]_{ij}LC \in \mathbb{M}^{r}$. Since $[\mathcal{A}]_{ij} \succeq 0$, to ensure $A-\sum_{v_j \in \mathcal{N}_{i}} [\mathcal{A}]_{ij}LC \in \mathbb{M}^{r}$, $\forall i\in \mathcal{I}$, it suffices to show that $A-l_{\mathrm{max}}LC \in \mathbb{M}^{ r}$. So, the positivity of the overall closed-loop system (6) is preserved if and only if $LC \succeq 0$, $BK \succeq 0$, $A - BK \in \mathbb{M}^{r}$ and $A-l_{\mathrm{max}}LC \in \mathbb{M}^{ r}$.

(ii) Consensus: To guarantee the consensus of the overall closed-loop system in (6), a well-known fact \cite{Wieland2011On} is that the consensus of system (6) is achievable if and only if $\tilde{A}-\lambda_{i}\tilde{B}L\tilde{C}$ is Hurwitz, $\forall i\in \mathcal{I} \backslash \{ 1\}$. By expanding $\tilde{A}-\lambda_{i}\tilde{B}L\tilde{C}$, we have
\begin{equation}\label{Amatrix3}
\tilde{A}-\lambda_{i}\tilde{B}L\tilde{C}=
\begin{bmatrix}
	A-BK      &    BK \\
	0                     &    A-\lambda_{i}LC
\end{bmatrix}.
\end{equation}
From (10), it is easy to see that $\tilde{A}-\lambda_{i}\tilde{B}L\tilde{C}$ is Hurwitz if and only if $A- \lambda_{i}LC \in \mathbb{H}^{r}$ and $A-BK \in \mathbb{H}^{r}$. Hence, the consensus of the multi-agent system in (\ref{SOF}) and (\ref{ABC}) is achievable if and only if $A-BK$ and $A- \lambda_{i}LC$, $ \forall i\in \mathcal{I} \backslash \{ 1\}$, are all Hurwitz.

The whole proof is completed.  $\hfill \hfill \square $
\vskip 0.2 cm
\noindent
\begin{myTheo}\label{CA4}
\textbf{Problem PCDMAS} is solvable if there exist a diagonal matrix $D>0$ , matrices $P>0$, $Q > 0$ and $S$ such that all the following statements hold:

\noindent
{\rm 1)} $PA^{\mathrm{T}} + AP - 2 {\rm Re}(\lambda_{2})PC^{\mathrm{T}}QCP < 0$,

\noindent
{\rm 2)} $A - l_{max}PC^{\mathrm{T}}QC \in \mathbb{M}^{r}$,

\noindent
{\rm 3)} $PC^{\mathrm{T}}QC \succeq 0$,

\noindent
{\rm 4)} $BS \succeq 0$,

\noindent
{\rm 5)} $AD - BS \in \mathbb{M}^{r}$,

\noindent
{\rm 6)} $AD - BS + DA^{\mathrm{T}} - S^{\mathrm{T}}B^{\mathrm{T}} < 0$.

\noindent
Under the conditions, $K = SD^{-1}$ and $L = PC^{\mathrm{T}}Q$.
\end{myTheo}

\textbf{Proof.}

According to the statement 1), we have $PA^{\mathrm{T}} + AP - 2{\rm Re}(\lambda_{2}) PC^{\mathrm{T}}QCP < 0$. Since $PC^{\mathrm{T}}QCP > 0$, the above inequality holds for any coefficient larger than $2{\rm Re}(\lambda_{2})$. Taking $L = PC^{\mathrm{T}}Q$, since $PA^{\mathrm{T}} + AP - 2{\rm Re}(\lambda_{i}) PC^{\mathrm{T}}QCP < 0$, $\forall i \in \mathcal{I}\backslash \{ 1\}$, we have $(A- \lambda_{i}LC)P + P (A - \lambda_{i}LC)^{\mathrm{*}} = PA^{\mathrm{T}} + AP - \lambda_{i} LCP - \lambda_{i}^{\mathrm{*}} PC^{\mathrm{T}}L^{\mathrm{T}} = PA^{\mathrm{T}} + AP - \lambda_{i} PC^{\mathrm{T}}QCP - \lambda_{i}^{\mathrm{*}} PC^{\mathrm{T}}QCP = PA^{\mathrm{T}} + AP - 2{\rm Re}(\lambda_{i}) PC^{\mathrm{T}}QCP < 0$, $\forall i \in \mathcal{I}\backslash \{ 1\}$. By the well-known result in \cite{Zhang2011Optimal}, we can conclude that $A - \lambda_{i}LC \in \mathbb{H}^{r}$, $\forall i \in \mathcal{I}  \backslash \{ 1\}$, which implies the statement 5) in \textbf{\rm Theorem 1}. Since $L = PC^{\mathrm{T}}Q$, the statements 2), 3) in \textbf{\rm Theorem 2} are equivalent to $A - l_{max}LC \in \mathbb{M}^{r}$ and $LC \succeq 0$, which implies the statements 3), 4) in \textbf{\rm Theorem 1}.

Taking $S = KD$, the statement 6) is equivalent to $AD - BKD + DA^{\mathrm{T}} - DK^{\mathrm{T}}B^{\mathrm{T}} < 0$. Hence we have $(A - BK)D + D(A - BK)^{\mathrm{T}} < 0$. By \textbf{\rm Lemma 2}, $A - BK$ is Hurwitz. From the statement 5), we have $AD - BKD \in \mathbb{M}^{r}$. As $D$ is a diagonal positive definite matrix, thus $A - BK \in \mathbb{M}^{r}$. So $A-BK$ is Hurwitz and Metzler, which implies the statement 2) in \textbf{\rm Theorem 1}. From the statement 4), we have $BKD \succeq 0$, so $BK \succeq 0$ as $D$ is a diagonal positive definite matrix. This implies the statement 1) in \textbf{\rm Theorem 1}. On the other hand, assuming that the statements 1), 2), 3) in \textbf{\rm Theorem 1} hold, by \textbf{\rm Lemma 2}, there must exist matrix $D > 0$ such that $(A - BK)D + D(A - KB)^{\mathrm{T}} < 0$. Taking $S = KD$, it is easy to see that, $AD - BS + DA^{\mathrm{T}} - S^{\mathrm{T}}B^{\mathrm{T}} < 0$, $BS \succeq 0$ and $AD - BS \in \mathbb{M}^{r}$.

The whole proof is completed.  $\hfill \hfill \square $
\vskip 0.2 cm
\begin{remark}
Notice that, \textbf{\rm Theorem 2} only needs us to focus on the Hurwitzness of $ A - \lambda_{2}LC$ instead of all the matrices $A - \lambda_{i}LC$, $\forall i \in \mathcal{I} \backslash \{ 1\} $. This fact is very useful in the later consensus design as it will greatly simplify the complexity of solving \textbf{Problem PCDMAS}.
\end{remark}
\vskip 0.2 cm
\begin{myTheo}\label{CD1}
\textbf{Problem PCDMAS} is solvable if there exist a diagonal matrix $D > 0$, matrices $P > 0$, $Q > 0$, $X > 0$, and $S$ such that all the following conditions hold:

\noindent
{\rm 1)} $PA^{\mathrm{T}} + AP - 2{\rm Re}(\lambda_{2})PC^{\mathrm{T}}QCX - 2{\rm Re}(\lambda_{2})XC^{\mathrm{T}}QCP + 2{\rm Re}(\lambda_{2})XC^{\mathrm{T}}QCX < 0$,

\noindent
{\rm 2)} $A - l_{max}PC^{\mathrm{T}}QC \in \mathbb{M}^{r}$,

\noindent
{\rm 3)} $PC^{\mathrm{T}}QC \succeq 0$,

\noindent
{\rm 4)} $BS \succeq 0$,

\noindent
{\rm 5)} $AD - BS \in \mathbb{M}^{r}$,

\noindent
{\rm 6)} $AD - BS + DA^{\mathrm{T}} - S^{\mathrm{T}}B^{\mathrm{T}} < 0$.

\noindent
Under the conditions, $K = SD^{-1}$ and $L = PC^{\mathrm{T}}Q$.
\end{myTheo}
\textbf{Proof.}

The proof of the statements 2), 3), 4), 5), 6) is similar to \textbf{\rm Theorem 2}. It suffices to show that, the statement 1) in \textbf{\rm Theorem 3} is equivalent to the statement 1) in \textbf{\rm Theorem 2}.

On one hand, if there exist matrices $P > 0$ and $X > 0$ such that $PA^{\mathrm{T}} + AP$ $- 2{\rm Re}(\lambda_{2})PC^{\mathrm{T}}QCX - 2{\rm Re}(\lambda_{2})XC^{\mathrm{T}}QCP + 2{\rm Re}(\lambda_{2})XC^{\mathrm{T}}QCX < 0$. Equivalently, we can obtain that, $PA^{\mathrm{T}} + AP - 2{\rm Re}(\lambda_{2})PC^{\mathrm{T}}QCP + 2{\rm Re}(\lambda_{2})(X - P)C^{\mathrm{T}}QC(X - P) < 0$. Because $(X - P)C^{\mathrm{T}}QC(X - P) \ge 0$, then we have that $PA^{\mathrm{T}} + AP - 2{\rm Re}(\lambda_{2})PC^{\mathrm{T}}QCP < 0$.
On the other hand, if $PA^{\mathrm{T}} + AP - 2{\rm Re}(\lambda_{2})PC^{\mathrm{T}}QCP < 0$ for some matrix $P > 0$, obviously there exists a matrix $X = P > 0$ such that $PA^{\mathrm{T}} + AP - 2{\rm Re}(\lambda_{2})PC^{\mathrm{T}}QCP + 2{\rm Re}(\lambda_{2})(X - P)C^{\mathrm{T}}QC(X - P) < 0$, i.e., $PA^{\mathrm{T}} + AP - 2{\rm Re}(\lambda_{2})PC^{\mathrm{T}}QCX - 2{\rm Re}(\lambda_{2})XC^{\mathrm{T}}QCP + 2{\rm Re}(\lambda_{2})XC^{\mathrm{T}}QCX < 0$.

The whole proof is completed. $\hfill \hfill \square $
\vskip 0.2 cm
\begin{remark}
In the above theorem, the matrix inequality in the statement 1) of Theorem 2 is linearized in an equivalent way by introducing the matrix $X$ and the matrix $X$ is assumed to be known. This enables us to use LMIs to solve \textbf{Problem PCDMAS}.
\end{remark}
\vskip 0.2 cm
Based on the results of \textbf{\rm Theorems 2} and \textbf{\rm 3}, a convex programming algorithm is developed as follows:
\vskip 0.2 cm
\newtheorem{algo}{Algorithm}
\begin{algo}
\
\begin{description}
    \item[Step 1:]\enskip Initialize $k=1$, $h = 1$, $Q^{(1)} = I$ and $\epsilon^{(0)}=0$. Obtain the initial matrix $P^{(1)} = U^{-1}$ by solving the following LMIs:
\noindent
\begin{equation*}
~~~~~~~~~~~
\begin{cases}
A^{\mathrm{T}}U + UA - 2{\rm Re}(\lambda_{2})C^{\mathrm{T}}Q^{(1)}C < 0\\
U > 0
\end{cases}~~~~~~~~~~~~~~~~
\end{equation*}
    \item[Step 2:]\enskip Set matrix $X^{(k)} = P^{(k)}$, minimize $\epsilon^{(h)}$
\begin{equation*}
\mathrm{s.t. }
\begin{cases}
PC^{\mathrm{T}}Q^{(k)}C \succeq 0\\
A - l_{max}PC^{\mathrm{T}}Q^{(k)}C \in \mathbb{M}^{r}\\
\Psi_{k}(P) < \epsilon^{(h)} I
\end{cases}
\end{equation*}
with respect to the matrix $P > 0$, where the matrix function $\Psi_{k}(P)$ is defined as $\Psi_{k}(P) := PA^{\mathrm{T}} + AP - 2{\rm Re}(\lambda_{2})PC^{\mathrm{T}}Q^{(k)}CX^{(k)} - 2{\rm Re}(\lambda_{2})X^{(k)}C^{\mathrm{T}}Q^{(k)}CP + 2{\rm Re}(\lambda_{2})X^{(k)}C^{\mathrm{T}}Q^{(k)}CX^{(k)}$.
    \item[Step 3:]\enskip If $\epsilon^{(h)} \preceq 0$, go to \textbf{\rm Step 8}. Otherwise, go to \textbf{\rm Step 4}.
    \item[Step 4:]\enskip If $|\epsilon^{(h)}-\epsilon^{(h-1)}|/\epsilon^{(h)} \prec \xi$, where $\xi$ is a prescribed tolerance, this algorithm fails to find the desired solution, \textbf{\rm STOP}. Otherwise, set $k=k+1$ and $h = h + 1$, update
 $P^{(k)} = P$, then go to \textbf{\rm Step 5}.
    \item[Step 5:]\enskip Minimize $\epsilon^{(h)}$
    \begin{equation*}
\mathrm{s.t. }
\begin{cases}
P^{(k)}C^{\mathrm{T}}QC \succeq 0\\
A - l_{max}P^{(k)}C^{\mathrm{T}}QC \in \mathbb{M}^{r}\\
\Lambda_{k}(Q) < \epsilon^{(h)} I
\end{cases}
\end{equation*}
with respect to the matrix $Q > 0$, where the matrix function $\Lambda_{k}(Q)$ is defined as $\Lambda_{k}(Q) := P^{(k)}A^{\mathrm{T}} + AP^{(k)} - 2{\rm Re}(\lambda_{2})P^{(k)}C^{\mathrm{T}}QCP^{(k)}$.
\item[Step 6:]\enskip If $\epsilon^{(h)}\preceq 0$, go to \textbf{\rm Step 8}. Otherwise, go to \textbf{\rm Step 7}.
\item[Step 7:]\enskip If $|\epsilon^{(h)}-\epsilon^{(h-1)}|/\epsilon^{(h)} \prec \xi$, where $\xi$ is a prescribed tolerance, then this algorithm fails to find the desired solution, \textbf{\rm STOP}. Otherwise, update
 $Q^{(k)} = Q$ and $h = h + 1$, then go to \textbf{\rm Step 2}.
    \item[Step 8:]\enskip Obtain matrices $D$ and $S$ by solving the following constraints:
\begin{equation*}
~~~~~~~~~~~~
\begin{cases}
BS \succeq 0\\
AD - BS \in \mathbb{M}^{r}\\
AD - BS + DA^{\mathrm{T}} - S^{\mathrm{T}}B^{\mathrm{T}} < 0
\end{cases}~~~~~~~~~~~~~~~~~~~~~~~
\end{equation*}
with respect to variables: diagonal matrix $D > 0$ and matrix $S$.
    \item[Step 9:]\enskip The feedback gain matrices $K$ and $L$ can be obtained as $K = SD^{-1}$ and $L = PC^{\mathrm{T}}Q$. \textbf{\rm STOP}.
\end{description}
\end{algo}
\vskip 0.2 cm
\begin{remark}
In Step 1 of the above algorithm, the matrix $P$ is initialized as $P^{(1)} = U^{-1}$ such that, $P^{(1)}A^{\mathrm{T}} + AP^{(1)} - 2{\rm Re(\lambda_{2})}P^{(1)}C^{\mathrm{T}}Q^{(1)}CP^{(1)} < 0$. Taking $L = P^{(1)}C^{\mathrm{T}}Q^{(1)}$, by \textbf{\rm Theorem 2}, we have that $A - \lambda_{i}LC$, $\forall i \in \mathcal{I} \backslash \{ 1\}$, are Hurwitz and thus such an initialization is reasonable. In \textbf{\rm Step 2}, $X^{(k)}$ is updated as $P^{(k)}$, $P^{(k+1)}$ always minimizes $\alpha(\Psi_{k}(P))$ and $Q^{(k)}$ minimizes $\alpha(\Lambda_{k}(Q))$, $\forall k \succeq 1$. Observing that, $\Lambda_{k+1}(Q^{k}) =  \Psi_{k}(P^{(k+1)})- 2{\rm Re}(\lambda_{2})(X^{(k)} - P^{(k+1)})C^{\mathrm{T}}Q^{(k)}C(X^{(k)} - P^{(k+1)})$, thus we have $ \alpha(\Lambda_{k+1}(Q^{k}))\preceq \alpha(\Psi_{k}(P^{k+1}))$ as $\Psi_{k}(P^{(k+1)}) - \Lambda_{k+1}(Q^{k}) \geq 0$. Moreover, because $\Psi_{k}(P^{(k)}) = \Lambda_{k}(Q^{k})$, we finally obtain that $ \alpha(\Lambda_{k+1}(Q^{k+1}))\preceq \alpha(\Lambda_{k+1}(Q^{k}))\preceq \alpha(\Psi_{k}(P^{(k+1)})) \preceq \alpha(\Psi_{k}(P^{(k)})) = \alpha(\Lambda_{k}(Q^{k}))$, $\forall k \succeq 1$. So we always have $\epsilon^{(h+1)} \preceq \epsilon^{(h)}$ for $h \succeq 1$ during the iterative process, which guarantees the convergence of \textbf{\rm Algorithm 1}. Meanwhile, positivity of the multi-agent system is also preserved since matrices $P$, $Q$ and $D$ only move in the feasible region during the iterations.
\end{remark}
\vskip 0.2 cm
\begin{remark}
Notice that, \textbf{\rm Algorithm 1} only needs to use eigenvalue $\lambda_{2}$ of the Laplacian matrix to design the feedback gain matrices $K$ and $L$, which provides an efficient way to solve \textbf{Problem PCDMAS}.
\end{remark}

\section{Numerical Simulation}
In this section, we use an example of directed multi-agent system to verify the effectiveness of the derived results and algorithm in this paper.
\begin{figure}[h]
	\centering
	\includegraphics[width=8cm]{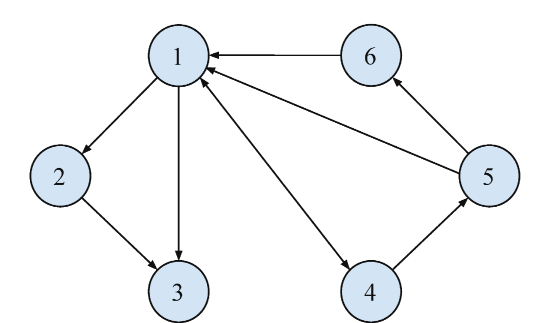}
	\caption{Directed communication graph}
	\label{fig:figure1}
\end{figure}
%\begin{figure}[!]
%	\centering
%	\includegraphics[width=8.8cm]{states.eps}
%	\caption{Phase plot of state $x_{i}$ in the multi-agent system.}
%	\label{fig:figure1}
%\end{figure}
\\Consider a multi-agent system in (2) with 6 agents and the following system matrices:
$$
A=
\begin{bmatrix}
	-3 & 2 & 3 \\
	1 & -4 & 2\\
    2 & 1 & -3
\end{bmatrix},~~~
B=
\begin{bmatrix}
	3 & 0       \\
	1 & 0 \\
    2 & 2
\end{bmatrix},~~~
C =
\begin{bmatrix}
	2 & 0 & 0\\
    0 & 2 & 3
\end{bmatrix}.
$$
The graph in \textbf{\rm Figure 1} is used to model the communication topology of the above multi-agent system. The associated Laplacian matrix of the multi-agent system is:
$$
\mathcal{L} =
\begin{bmatrix}
	3&0&0&-1&-1&-1\\
    -1&1&0&0&0&0\\
    -1&-1&2&0&0&0\\
    -1&0&0&1&0&0\\
    0&0&0&-1&1&0\\
    0&0&0&0&-1&1
\end{bmatrix}
.
$$
Using the LMI Toolbox of MATLAB, the \textbf{\rm Algorithm 1} is implemented and the following results are obtained. \\
The matrix $P$ is initialized as $P^{(1)} =$
$$
\begin{bmatrix}
        0.60314 &    -0.02879&   -0.0068583\\
     -0.02879    &  0.33887   &   0.03418\\
   -0.0068583     & 0.03418    &  0.60242
\end{bmatrix}
.
$$
The matrix $Q$ is initialized as $Q^{(1)} = I$. Obviously, the pair $(P^{(1)}, Q^{(1)})$ is not feasible since $P^{(1)}C^{\mathrm{T}}Q^{(1)}C =$
$$
\begin{bmatrix}
	2.4126  &   -0.15631  &   -0.23446\\
     -0.11516&       1.5606&       2.3408\\
    -0.027433 &      3.7512 &      5.6269
\end{bmatrix}
 \nsucceq 0
$$
and $A - l_{max}P^{(1)}C^{\mathrm{T}}Q^{(1)}C =$
$$
\begin{bmatrix}
	 -10.238  &     2.4689   &    3.7034 \\
       1.3455  &    -8.6817   &   -5.0225 \\
       2.0823   &   -10.254    &  -19.881
\end{bmatrix}
 \notin \mathbb{M}^{r}.
$$
After several iterations, $\epsilon = -0.12567 \preceq 0$ and meanwhile we have $P =$
$$
\begin{bmatrix}
1.8327  & 5.0754\mathrm{e}{-06} &  2.4062\mathrm{e}{-06}\\
5.0754\mathrm{e}{-06}    &   5.5882    &  -3.6514\\
2.4062\mathrm{e}{-06}    &  -3.6514    &   2.4898
\end{bmatrix}
$$
and $Q = $
$$
\begin{bmatrix}
        9.9978  &   0.090679\\
       0.090679  &     0.9999
\end{bmatrix}.
$$
Observe that, $PC^{\mathrm{T}}QC =$
$$
\begin{bmatrix}
      73.292   &   0.66478   &   0.99717 \\
      0.040505  &     0.4444  &     0.6666 \\
      0.030322   &    0.3333   &   0.49995
\end{bmatrix}
 \succeq 0
$$
and $A - l_{max}PC^{\mathrm{T}}QC =$
$$
\begin{bmatrix}
	 -222.87 &   0.0056585  &  0.0084877 \\
      0.87848 &     -5.3332  &  0.00020407 \\
        1.909  & 0.00010217   &   -4.4998
\end{bmatrix}
  \in \mathbb{M}^{r}.
$$
So, the obtained matrix pair $(P,~Q)$ is feasible. \\
Then the feedback gain matrices $K$ and $L$ are obtained as:
$$
K =
\begin{bmatrix}
	0.30383   &   0.36377  &    0.59561\\
      0.20697  &   -0.18142 &    -0.14544
\end{bmatrix}
$$
and
$$
L =
\begin{bmatrix}
         36.646  &    0.33239 \\
         0.020252 &      0.2222 \\
         0.015161  &    0.16665
\end{bmatrix}.
$$
To better illustrate the effectiveness of our algorithm, the states of the six augmented systems in the above example are respectively initialized as,
\begin{equation*}\label{}
{\tilde{x}}_{1}(0) =
\begin{bmatrix}
20\\
1\\
1\\
15\\
0.5\\
0.5
\end{bmatrix},~~~
{\tilde{x}}_{2}(0) =
\begin{bmatrix}
15\\
2\\
2\\
10\\
1\\
1
\end{bmatrix},~~~
{\tilde{x}}_{3}(0) =
\begin{bmatrix}
13\\
5\\
1\\
10\\
2\\
0.5
\end{bmatrix},
\end{equation*}
and
\begin{equation*}
{\tilde{x}}_{4}(0) =
\begin{bmatrix}
14 \\
10\\
5 \\
10\\
6\\
3
\end{bmatrix},~~~
{\tilde{x}}_{5}(0) =
\begin{bmatrix}
12\\
8\\
6\\
8\\
4\\
2
\end{bmatrix},~~~
{\tilde{x}}_{6}(0) =
\begin{bmatrix}
9\\
6\\
4\\
7\\
1\\
0.5
\end{bmatrix}.
\end{equation*}

%$
%[{\tilde{x}}_{1}(0), ~{\tilde{x}}_{2}(0), ~{\tilde{x}}_{3}(0), ~{\tilde{x}}_{4}(0), ~{\tilde{x}}_{5}(0), ~{\tilde{x}}_{6}(0)]^{\mathrm{T}} =
%$
%$$
%\begin{bmatrix}
%20 &1 &1& 15 &0.5 &0.5\\
%15 &2 &2 &10 &1 &1\\
%13 &5 &1 &10 &2  &0.5\\
%14 &10& 5 &10& 6& 3\\
%12 &8& 6 &8 &4 &2\\
%9 &6& 4 &7 &1 &0.5
%\end{bmatrix}.
%$$
With the obtained feedback gain matrices $K$ and $L$, the multi-agent system finally reaches consensus and the phase plots are shown in \textbf{\rm Figures 2} and \textbf{\rm 3}.
\begin{figure}[h]
	\centering
\psfrag{x1}{$x_{i1}(t)$}
\psfrag{x2}{$x_{i2}(t)$}
\psfrag{x3}{$x_{i3}(t)$}
	\includegraphics[width=8.3cm]{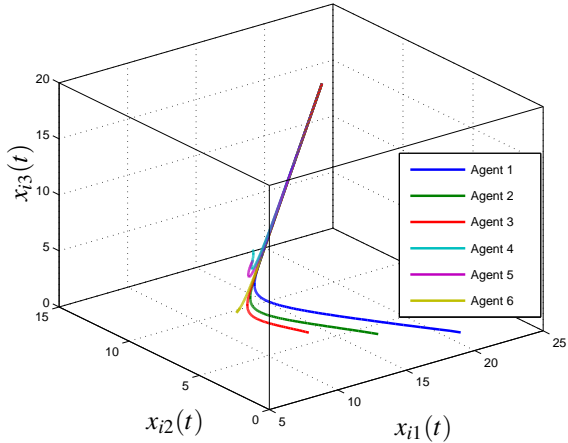}
	\caption{Phase plot of the state $x_{i}(t)$ in the multi-agent system}
	\label{fig:figure1}
\end{figure}
\begin{figure}[h]
	\centering
\psfrag{e1}{$e_{i1}(t)$}
\psfrag{e2}{$e_{i2}(t)$}
\psfrag{e3}{$e_{i3}(t)$}
	\includegraphics[width=8.3cm]{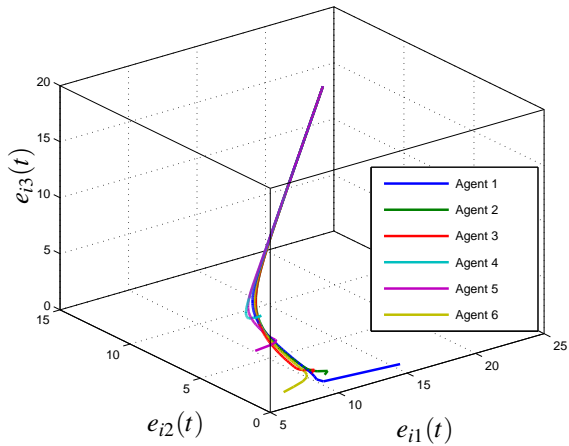}
	\caption{Phase plot of the feedback signal $e_{i}(t)$ in the multi-agent system}
	\label{fig:figure1}
\end{figure}
\section{Conclusion}
This paper has studied the positive consensus problem of directed multi-agent systems. Based on the results in positive systems theory and the consensus issue, the positive consensus of directed multi-agent systems has been analyzed. Some new results have been derived in the form of LMIs and a convex programming algorithm has been developed to design appropriate observer-type protocols such that the multi-agent system is able to reach consensus with its state trajectory always remaining in the non-negative orthant. Finally, the simulations have illustrated the effectiveness of the derived results and algorithm.
\vspace{8.8pt}
\bibliographystyle{abbrv}
\bibliography{reference}
%%%%%%%%%%%%%%%%%%%%%%%%%%%%%%%%%%%%%%%%%%%%%%%%%%%%%%%%%%%%%%%%%%%%%%%%%%%%%%%%

\end{document}